# Observability and Controllability of Second Order Linear Time Invariant Systems and Kalman Type Conditions

**Elimhan N. Mahmudov**[a,b]

[a]*Department of Mathematics, Istanbul Technical University, Istanbul, Turkey,*

[b]*Azerbaijan National Academy of Sciences Institute of Control Systems, Baku,*

*Azerbaijan.*

e-mail: *elimhan22@yahoo.com*

**Abstract.** In the present paper we consider controllability and observability of second order linear time invariant systems in matrix form. Without reducing into first order systems we show how the classical conditions for first order linear systems can be generalized to this case. In term of Kalman type criterions these concepts are investigated for second order discrete and continuous time linear systems. It should be pointed out that by repeated differentiation of state and output vector-functions we derive two different systems of linear algebraic equations. Then the initial values $x_0, x_1$ and input functions can be determined uniquely from these systems if and only if the observability and controllability matrices have full rank, respectively. Also the transfer function of the second order continuous-time linear state-space system is constructed. A numerical example is given to illustrate the feasibility and effectiveness of the theoretic results obtained.

1. Introduction

The behaviour of many important systems in practical engineering practice is governed by second order controllable systems, our objective in the present paper is to contribute to the analysis of the controllability and observability of linear discrete/continuous time of second order systems; second order linear systems theory has played an important role in many technology advancements, for example, in aerospace, communications, automotive, and computer engineering, forced oscillations and resonance phenomenon, earthquake–induced vibrations of multi-storey buildings, etc. We recall that controllability and observability play a crucial role in many control problems, involving time invariant and time-varying systems and are dual aspects of the same problem and represent two major concepts of modern control system theory (Zadeh, L.A. and Desoer, C.A. 1963). These concepts were introduced by Kalman (Kalman, R.E. 1960). They can be roughly defined as follows; in order to be able to do whatever we want with the given dynamic system and to see what is going on inside under control input, the system must be controllable and observable. In this paper we show that the concepts of controllability and observability are related to linear systems of algebraic equations. It is well known that a solvable system of linear algebraic equations has a solution if and only if the rank of the system matrix is full. Observability and controllability tests will be connected to the rank tests of certain matrices, i.e. the controllability and observability matrices.

The works (Fucheng Liao, Yu Wang, Yanrong Lu and Jiamei Deng, 2017; Gałkowski, K., Sulikowski, B., Paszke, W., *et al.* 2001; Gen, Qi Xu, Chao Liu and Siu Pang Yung 2008; Isidori, A. and Ruberti, A. 1976; Lin, L., Stefanescu, A., Su, R., Wang, W. and Shehabinia, A. 2014; Lu, J., Zhong, J., Daniel, Ho, *et.al.* 2016; Mehmet Önder Efe, Okyay Kaynak and Bogdan M. 2000; Memet Kule, 2015; Silverman, L.M. and Meadows, H.E. 1967; Tang, Y., Gao, H. Yang Shi, *et al.* 2015; Weili Yan, Chunling Du and Chee Khiang Pang, 2016; Yanling Wei, Jianbin Qiu, Peng Shi and Chadli, M. 2016; Zhendong Sun and Shuzhi Sam Ge, 2003; Zhijian Ji, Hai Lin, Gang Feng and Xiaoxia Guo, 2010) focuses on different first order time-invariant and time-varying finite-dimensional systems, covering both continuous and discrete-time topics.

In the monograph and papers (Mahmudov, E.N. 2011; Mahmudov, E.N. 2015a;

Mahmudov, E.N. 2015b; Mahmudov, E.N. 2018a; Mahmudov, E.N. 2018b; Mahmudov E.N. 2019) of author different problems of optimal control theory with fixed and varying time interval given by convex, nonconvex higher order discrete and differential inclusions involving second order linear systems are considered. Necessary and sufficient optimality conditions under the transversality condition are proved incorporating the Euler-Lagrange and Hamiltonian type inclusions.

It is well known that performing a classical transformation to first order form may destroy some controllability and observability properties. To avoid this, we use separation of the equations involving even and odd-orders finite differences and derive a second order discrete-time/continuous time system in the state space that do not destroy the controllability and observability properties.

The remainder of this paper is organized as follows: In Section 2 we consider a linear, time invariant, second order discrete-time system in the state space form with output measurements. The method applied in the present analysis is not standard, and its exposition in the framework of investigated problem is new. It is shown that using even and odd-orders finite differences and observability matrix the initial values $x_0, x_1$ can be determined. In Section 3 for the purpose of studying observability of second order continuous system, we consider an input-free system with the corresponding measurements and show how the classical observability conditions for first order systems can be generalized to this case. Repeatedly differentiating of output vector-function we derive two different system of linear algebraic equations. Then the initial values $x_0, x_1$ can be determined uniquely from these systems if and only if the observability matrix has full rank. Section 4 deals with the controllability of second order time-invariant linear discrete systems; it is possible to find a control sequence such that the given point can be reached from any initial state in a finite time. The system $x_{t+2} = \tilde{A}_0 x_t + \tilde{A}_1 x_{t+1} + \tilde{B} u_t$ is controllable if and only if $\operatorname{rank} \mathcal{C}(\tilde{A}_0, \tilde{A}_1, \tilde{B}) = n$. Section 5 devoted to study the evolution of the concept of controllability of second order linear continuous systems. We show how to find a control input that will transfer our system from any initial state to any final state. For a vector input system $x'' = A_0 x + A_1 x' + Bu$ the above discussion produces the same relation with the controllability matrix and with the input vector $u(t)$. We remarked that the observability results of Section 3 can be obtained from duality relation between observability and controllability. It is shown that for controllability of this system $\operatorname{rank} \mathcal{C}(A_0, A_1, B) = n$. Besides, if $A_1$ is a zero matrix, then $\mathcal{C}(A_0, O_{n \times n}, B) = \mathcal{C}(A_0, B)$. Numerical examples are presented to illustrate the theoretical results; we present some examples corresponding to different cases of checking the controllability and observability conditions. In particular, transfer function $H(s)$ is constructed and in an example involving single-input, single-output system, by zero-pole cancellations principle of the transfer function we show that system is observable and controllable.

## 2. Observability of Second Order Discrete Systems

This section is devoted to observability of controllable second order discrete systems. Criteria for determining observability for these systems are developed. Thus, let us consider a linear, time invariant, second order discrete-time system in the state space form

$$x_{t+2} = \tilde{A}_0 x_t + \tilde{A}_1 x_{t+1} + \tilde{B} u_t, \ x_0 = a_0, \ x_1 = a_1, \tag{1}$$

with output measurements

$$y_t = \tilde{C} x_t, \ t = 0, 1, 2, .... \tag{2}$$

where $x_t = (x_t^1, ..., x_t^n)^T \in \mathbb{R}^n$, $y_t = (y_t^1, ..., y_t^p)^T \in \mathbb{R}^p$, $a_0, a_1$ unknown vectors, $\tilde{A}, \tilde{B}, \tilde{C}$ are constant matrices of appropriate dimensions; $\tilde{A} \in \mathbb{R}^{n,n}$, $\tilde{B} \in \mathbb{R}^{n,r}$, $\tilde{C} \in \mathbb{R}^{p,n}$. Here $\mathbb{R}^{n,\ell}$ denotes the vector space of $n \times \ell$ real matrices, $x_t$ is the state, $u_t$ is the unknown input vector, and $y_t$ is the output of the system. Now using only information from the output measurements (2) we will solve the observability problem for the state space variables defined in (1). We say that the second order linear discrete-time system modelled by (1) and (2) *is observable* at time $t = 0$ if there exits some $t_1 > 0$ such that the state $(x_0, x_1) = (a_0, a_1)$ at time $t = 0$ can be uniquely determined from the knowledge of $u_t, y_t, t = 0, 1, ..., t_1$.

For a known $a_0, a_1$, of course, the recursion (1) gives us complete knowledge about the state variables at any discrete-time instant. Thus, our main problem is to determine from the state measurements the initial state vector $x_0 = a_0$, $x_1 = a_1$. The developments applied in this section are not standard, and their exposition in the framework of problem (1) is new.

Obviously, the $n$-dimensional vector has unknown components, and we try to determine an initial values $a_0$ and $a_1$ using output measurements $y_t$. Take $t = 0,1,2,...$ in (1) and (2), and generate the following sequence

$$y_0 = \tilde{C}x_0,\ y_1 = \tilde{C}x_1;\ \ y_2 = \tilde{C}\tilde{A}_0 x_0 + \tilde{C}\tilde{A}_1 x_1 + \tilde{C}\tilde{B}u_0,\ y_3 = \tilde{C}\tilde{A}_1\tilde{A}_0 x_0 + \tilde{C}\left(\tilde{A}_0 + \tilde{A}_1^2\right)x_1 + \tilde{C}\left(\tilde{A}_1\tilde{B}u_0 + \tilde{B}u_1\right),$$
$$y_4 = \tilde{C}(\tilde{A}_0^2 + \tilde{A}_1^2\tilde{A}_0)x_0 + \tilde{C}(\tilde{A}_1\tilde{A}_0 + \tilde{A}_0\tilde{A}_1 + \tilde{A}_1^3)x_1 + \tilde{C}\left[(\tilde{A}_0 + \tilde{A}_1^2)\tilde{B}u_0 + \tilde{A}_1\tilde{B}u_1 + \tilde{B}u_2\right],\ \ldots .$$

(3)

Now using recurrence relation below let us rewrite (3) in more suitable form. Take $S_0 = \tilde{A}_0$, $P_0 = \tilde{A}_1$ and denote $S_k = P_{k-1}\tilde{A}_0$; $P_k = S_{k-1} + P_{k-1}\tilde{A}_1\ (k = 1,2,3,...)$. Then it is not hard to see that the coefficient matrices of $x_0$ and $x_1$ in (3) can be expressed in terms of two couple of matrices $S_k$ and $P_k$ defined as follows

$$\begin{aligned} S_1 &= \tilde{A}_1\tilde{A}_0;\ S_2 = \tilde{A}_0^2 + \tilde{A}_1^2\tilde{A}_0, & P_1 &= \tilde{A}_0 + \tilde{A}_1^2;\ P_2 = \tilde{A}_1\tilde{A}_0 + \tilde{A}_0\tilde{A}_1 + \tilde{A}_1^3, \\ S_3 &= \tilde{A}_1\tilde{A}_0^2 + \tilde{A}_0\tilde{A}_1\tilde{A}_0 + \tilde{A}_1^3\tilde{A}_0, & P_3 &= \tilde{A}_0^2 + \tilde{A}_1^2\tilde{A}_0 + \tilde{A}_1\tilde{A}_0\tilde{A}_1 + \tilde{A}_0\tilde{A}_1^2 + \tilde{A}_1^4, \\ &\ldots\ldots\ldots\ldots\ldots\ldots\ldots & &\ldots\ldots\ldots\ldots\ldots\ldots\ldots\ldots \end{aligned}$$ (4)

On the other hand denoting $M_k = \tilde{A}_0 M_{k-2} + \tilde{A}_1 M_{k-1}\ (M_0 = \tilde{B}, M_1 = \tilde{A}_1\tilde{B})$, $k = 2,3,4,...$ it can be easily checked that

$$M_2 = \tilde{A}_0 M_0 + \tilde{A}_1 M_1 = (\tilde{A}_0 + \tilde{A}_1^2)\tilde{B}; \ M_3 = \tilde{A}_0 M_1 + \tilde{A}_1 M_2 = (\tilde{A}_0\tilde{A}_1 + \tilde{A}_1\tilde{A}_0 + \tilde{A}_1^3)\tilde{B},$$
$$M_4 = \tilde{A}_0 M_2 + \tilde{A}_1 M_3 = (\tilde{A}_0^2 + \tilde{A}_0\tilde{A}_1^2 + \tilde{A}_1\tilde{A}_0\tilde{A}_1 + \tilde{A}_1^2\tilde{A}_0 + \tilde{A}_1^4)\tilde{B}, \cdots \cdots \quad (5)$$

Then using mathematical induction in terms of $S_k, P_k, M_k$ the sequence of equations (3) has the following remarkable representation

$$y_0 = \tilde{C}x_0; \ y_1 = \tilde{C}x_1; \ y_2 = \tilde{C}S_0 x_0 + \tilde{C}P_0 x_1 + \tilde{C}M_0 u_0; \ y_3 = \tilde{C}S_1 x_0 + \tilde{C}P_1 x_1 + \tilde{C}(M_1 u_0 + M_0 u_1),$$
$$\cdots, \ y_{2n-1} = \tilde{C}S_{2n-3} x_0 + \tilde{C}P_{2n-3} x_1 + \tilde{C}(M_{2n-3} u_0 + M_{2n-2} u_1 + \cdots + M_0 u_{2n-3}).$$
$$(6)$$

Now we rewrite the equation (6) in the following relevant matrix form

$$\begin{bmatrix} y_0 \\ y_1 \\ y_2 - \tilde{C}M_0 u_0 \\ \vdots \\ y_{2n-1} - \tilde{C}(M_{2n-3} u_0 + M_{2n-2} u_1 + \cdots + M_0 u_{2n-3}) \end{bmatrix}_{(2np) \times 1} = \begin{bmatrix} \tilde{C} & O_{p \times n} \\ O_{p \times n} & \tilde{C} \\ \tilde{C}S_0 & \tilde{C}P_0 \\ \vdots & \vdots \\ \tilde{C}S_{2n-3} & \tilde{C}P_{2n-3} \end{bmatrix} \cdot \begin{bmatrix} x_0 \\ x_1 \end{bmatrix}, \quad (7)$$

where $O_{\ell \times m}$ is $\ell \times m$ zero matrix. It is well known from linear algebra that the system of linear algebraic equations with $2n$ unknowns, (6) has a unique solutions $x_0$ and $x_1$ if and only if the system matrix has rank $2n$.

$$\operatorname{rank} \begin{bmatrix} \tilde{C} & O_{p \times n} \\ O_{p \times n} & \tilde{C} \\ \tilde{C}S_0 & \tilde{C}P_0 \\ \vdots & \vdots \\ \tilde{C}S_{2n-3} & \tilde{C}P_{2n-3} \end{bmatrix}_{(2np) \times 2n} = \operatorname{rank} \mathcal{O}(\tilde{A}_0, \tilde{A}_1, \tilde{C}) = 2n. \quad (8)$$

Thus, the initial values $x_0, x_1$ can be determined if the so-called observability matrix $\mathcal{O}(\tilde{A}_0, \tilde{A}_1, \tilde{C})$ has rank $2n$. Thus we have the following theorem.

**Theorem 2.1** The second order linear discrete-time system (1) with measurements (2) is observable if and only if the observability matrix (8) has rank equal to $2n$.

Let us outline the following particular cases. First, suppose that in the second order linear discrete-time system (1) $\tilde{A}_1$ is the null matrix.

**Corollary 2.1** The second order linear discrete-time system (1) with the null matrix $\tilde{A}_1$ and with the measurements (2) is observable if and only if $\operatorname{rank} \mathcal{O}(\tilde{A}_0, O_{n \times n}, \tilde{C}) = n$.

*Proof* Since $P_k = S_{k-1} + P_{k-1}\tilde{A}_1$, $S_k = P_{k-1}\tilde{A}_0 \ (k=1,2,3,...)$ it follows that $P_k = S_{k-1}$ and then $S_k = S_{k-2}\tilde{A}_0$. Thus taking into account that $S_0 = \tilde{A}_0$, $S_1 = O_{n \times n}$ for all odd numbers $2k-1$ we have $S_{2k-1} = O_{n \times n}$, $k=1,2,...,n$. Therefore for an even indices $n = 2k$ it follows that $P_{2k} = S_{2k-1} = O_{n \times n}$. Hence an observability matrix (8) has the following form

$$\mathcal{O}(\tilde{A}_0, \tilde{A}_1, C) = \mathcal{O}(\tilde{A}_0, O_{n \times n}, \tilde{C}) = \begin{bmatrix} \tilde{C} & O_{p \times n} \\ O_{p \times n} & \tilde{C} \\ \tilde{C}\tilde{A}_0 & O_{p \times n} \\ \vdots & \vdots \\ \tilde{C}\tilde{A}_0^{n-1} & O_{p \times n} \\ O_{p \times n} & \tilde{C}\tilde{A}_0^{n-1} \end{bmatrix}.$$

On the other hand in view of $M_k = \tilde{A}_0 M_{k-2} + \tilde{A}_1 M_{k-1} \ (M_0 = \tilde{B}, M_1 = O_{n \times r}), (k=2,3,4,...)$ we have $M_k = \tilde{A}_0 M_{k-2}$ and the system (6) can be transformed into the system

$$y_0 = \tilde{C} x_0;\ y_1 = \tilde{C} x_1;\ y_2 = \tilde{C}\tilde{A}_0 x_0 + \tilde{C}\tilde{B} u_0,\ y_3 = \tilde{C}\tilde{A}_0 x_1 + \tilde{C}\tilde{B} u_1;\ y_4 = \tilde{C}\tilde{A}_0^2 x_0 + \tilde{C}(\tilde{A}_0 \tilde{B} u_0 + \tilde{B} u_2),$$
$$\cdots,\ y_{2n-2} = \tilde{C}\tilde{A}_0^{n-1} x_0 + \tilde{C}(\tilde{A}_0^{n-2}\tilde{B} u_0 + \tilde{A}_0^{n-3}\tilde{B} u_2 + \cdots + \tilde{B} u_{2n-4}),$$

$$y_{2n-1} = \tilde{C}\tilde{A}_0^{n-1} x_1 + \tilde{C}(\tilde{A}_0^{n-2}\tilde{B} u_1 + \tilde{A}_0^{n-3}\tilde{B} u_3 + \cdots + \tilde{B} u_{2n-3}). \tag{9}$$

Now we can separate the equations (9) involving even and odd-orders finite differences in the following relevant matrix form

$$\begin{bmatrix} y_0 \\ y_2 - \tilde{C}\tilde{B} u_0 \\ \vdots \\ y_{2(n-1)} - \tilde{C}\left(\tilde{A}_0^{n-2}\tilde{B} u_0 + \tilde{A}_0^{n-3}\tilde{B} u_2 + \cdots + \tilde{B} u_{2n-4}\right) \end{bmatrix}_{(np)\times 1} = \begin{bmatrix} \tilde{C} \\ \tilde{C}\tilde{A}_0 \\ \vdots \\ \tilde{C}\tilde{A}_0^{n-1} \end{bmatrix}_{(np)\times n} x_0;$$

$$\begin{bmatrix} y_1 \\ y_3 - \tilde{C}\tilde{B} u_1 \\ \vdots \\ y_{2n-1} - \tilde{C}\left(\tilde{A}_0^{n-2}\tilde{B} u_1 + \tilde{A}_0^{n-3}\tilde{B} u_3 + \cdots + \tilde{B} u_{2n-3}\right) \end{bmatrix}_{(np)\times 1} = \begin{bmatrix} \tilde{C} \\ \tilde{C}\tilde{A}_0 \\ \vdots \\ \tilde{C}\tilde{A}_0^{n-1} \end{bmatrix}_{(np)\times n} x_1.$$

Then the above conditions converted to the Kalman criterion $\text{rank}\,\mathcal{O}(\tilde{A}_0, \tilde{C}) = n$, where

$$\mathcal{O}(\tilde{A}_0, O_{n\times n}, \tilde{C}) = \begin{bmatrix} \tilde{C} \\ \tilde{C}\tilde{A}_0 \\ \vdots \\ \tilde{C}\tilde{A}_0^{n-1} \end{bmatrix}_{(np)\times n}.$$

**Corollary 2.2** The second order linear discrete-time system (1) with zero matrix $\tilde{A}_0$ and with measurements (2) is observable if and only if the observability matrix (9) has rank equal to $n$, i.e. $\text{rank}\,\mathcal{O}(\tilde{A}_1, O_{n\times n}, \tilde{C}) = n$.

*Proof* In this case it can be easily seen that $S_k = P_{k-1}\tilde{A}_0 = O_{n\times n}$ for all $k = 1, 2, \ldots, 2n-3$

and $P_k = P_{k-1}\tilde{A}_1 = P_{k-2}\tilde{A}_1^2 = \ldots = P_0\tilde{A}_1^k = \tilde{A}_1^{k+1}$. Then taking higher order derivatives until

$n$ in the system (6) we have

$$y_1 = \tilde{C}x_1;\ y_2 = \tilde{C}\tilde{A}_1 x_1 + \tilde{C}\tilde{B}u_0,\ y_3 = \tilde{C}\tilde{A}_1^2 x_1 + \tilde{C}\left(\tilde{A}_1\tilde{B}u_0 + \tilde{B}u_1\right),\ \cdots\cdots, \qquad (10)$$
$$y_n = \tilde{C}\tilde{A}_1^{n-1} x_1 + \tilde{C}\left(\tilde{A}_1^{n-2}\tilde{B}u_0 + \tilde{A}_1^{n-3}\tilde{B}u_1 + \cdots + \tilde{B}u_{n-2}\right).$$

In turn from (10) we conclude that the observability matrix (8) consists of the following

$$\mathcal{O}(\tilde{A}_1, \tilde{C}) = \begin{bmatrix} \tilde{C} \\ \tilde{C}\tilde{A}_1 \\ \vdots \\ \tilde{C}\tilde{A}_1^{n-1} \end{bmatrix}_{(np)\times n}.$$

Obviously, if $\text{rank}\,\mathcal{O}(\tilde{A}_1, \tilde{C}) = n$ we can define unique initial value $x_1$ of the system (12).

We note that this occurs because by substituting $\tilde{A}_0 = O_{n\times n}$ into system (1) we have a first order discrete-time equations instead of the second order.

**Example 2.1** Consider the following second order linear discrete-time system with

measurements

$$\begin{bmatrix} x_{t+2}^1 \\ x_{t+2}^2 \end{bmatrix} = \begin{bmatrix} 1 & 0 \\ 1 & -1 \end{bmatrix} \begin{bmatrix} x_t^1 \\ x_t^2 \end{bmatrix} + \begin{bmatrix} 0 & 1 \\ 1 & 2 \end{bmatrix} \begin{bmatrix} x_{t+1}^1 \\ x_{t+1}^2 \end{bmatrix},$$

$$y_t = \begin{bmatrix} 2 & 1 \end{bmatrix} \begin{bmatrix} x_t^1 \\ x_t^2 \end{bmatrix}; \tilde{A}_0 = \begin{bmatrix} 1 & 0 \\ 1 & -1 \end{bmatrix}, \tilde{A}_1 = \begin{bmatrix} 0 & 1 \\ 1 & 2 \end{bmatrix}, \tilde{C} = \begin{bmatrix} 2 & 1 \end{bmatrix}, \tilde{B} = O_{2\times 1}.$$

In this case $M_0 = O_{2\times 1}, M_1 = O_{2\times 1}$ and the equations (7) is simplified as follows

$$\begin{bmatrix} y_0 \\ y_1 \\ y_2 \\ y_3 \end{bmatrix} = \begin{bmatrix} \tilde{C} & O_{1\times 2} \\ O_{1\times 2} & \tilde{C} \\ \tilde{C}S_0 & \tilde{C}P_0 \\ \tilde{C}S_1 & \tilde{C}P_1 \end{bmatrix} \cdot \begin{bmatrix} x_0^1 \\ x_0^2 \\ x_1^1 \\ x_1^2 \end{bmatrix},$$

where $S_0 = \tilde{A}_0$, $P_0 = \tilde{A}_1$, $S_1 = \tilde{A}_1\tilde{A}_0$, $P_1 = \tilde{A}_0 + \tilde{A}_1^2$. Therefore $\tilde{C}S_0 = \begin{bmatrix} 3 & -1 \end{bmatrix}$, $\tilde{C}P_0 = \begin{bmatrix} 1 & 4 \end{bmatrix}$.

Besides

$$S_1 = \begin{bmatrix} 1 & -1 \\ 3 & -2 \end{bmatrix}, P_1 = \begin{bmatrix} 2 & 2 \\ 3 & 4 \end{bmatrix}, \tilde{C}S_1 = \begin{bmatrix} 5 & -4 \end{bmatrix}, \tilde{C}P_1 = \begin{bmatrix} 7 & 8 \end{bmatrix}, \quad \mathcal{O}(\tilde{A}_0, \tilde{A}_1, \tilde{C}) = \begin{bmatrix} 2 & 3 & 0 & 0 \\ 0 & 0 & 2 & 3 \\ 3 & -1 & 1 & 4 \\ 5 & -4 & 7 & 8 \end{bmatrix}.$$

Then using the cofactor expansion of $\mathcal{O}(\tilde{A}_0, \tilde{A}_1, \tilde{C})$ along the first or second row it is easy to compute that the determinant $\det \mathcal{O}(\tilde{A}_0, \tilde{A}_1, \tilde{C})$ of this $4\times 4$ square matrix is nonzero ($\det \mathcal{O}(A,C) = -175 \neq 0$) and the rows of the matrix are linearly independent. Thus $\operatorname{rank} \mathcal{O}(\tilde{A}_0, \tilde{A}_1, \tilde{C}) = 2n = 4$ i.e. the system under consideration is observable.

**Example 2.2** Suppose now that the corresponding system and measurement matrices of our problem are given as follows

$$\tilde{A}_0 = \begin{bmatrix} 3 & 2 \\ -2 & -1 \end{bmatrix}, \tilde{A}_1 = O_{2\times n}, C = \begin{bmatrix} 3 & 1 \end{bmatrix}. \text{ Obviously, } \mathcal{O}(\tilde{A}_0, O_{2\times 2}, \tilde{C}) = \begin{bmatrix} \tilde{C} \\ \tilde{C}\tilde{A}_0 \end{bmatrix} = \begin{bmatrix} 3 & 1 \\ -3 & -1 \end{bmatrix}.$$

Therefore, $\text{rank}\,\mathcal{O}(\tilde{A}_0, O_{2\times 2}, \tilde{C}) = 1$. It means that $\text{rank}\,\mathcal{O}(\tilde{A}_0, O_{2\times 2}, \tilde{C}) = 1 < n$, where $n = 2$ and the system is unobservable.

## 3. Observability of Second Order Continuous Systems

In this section we study an observability problem for a second order continuous-time linear state-space system

$$x''(t) = A_0 x(t) + A_1 x'(t) + Bu(t),\ x(0) = x_0,\ x'(0) = x_1 \tag{11}$$

with the corresponding measurements

$$y(t) = Cx(t) \tag{12}$$

of dimensions $x(t) \in \mathbb{R}^n$, $y(t) \in \mathbb{R}^p$ and matrices $A_0, A_1 \in \mathbb{R}^{n,n}$, $B \in \mathbb{R}^{n,r}$ and $C \in \mathbb{R}^{p,n}$, $x_0, x_1$ unknown vectors. Here $x(t) \in \mathbb{R}^n$, $u(t) \in \mathbb{R}^r$ and $y(t) \in \mathbb{R}^p$ are respectively the state vector, the unknown input vector (or input signal) and the output vector, $x(\cdot)$ and $u(\cdot)$ are continuously differentiable functions of orders $2n-1$ and $2n-3$, respectively. We shall conclude that the knowledge of $x_0, x_1$ is sufficient to determine at any time instant. This difficulty is also related to the one we encountered in the previous section when dealing with the second order time invariant discrete systems. The problem is to find $x_0, x_1$ from the available measurements (12) that is, system (11) is observable on an interval $[0, t_1]$ if any initial state $(x_0, x_1)$ at $t = 0$ can be determined from knowledge of the system output $y(t)$ and input $u(t)$ over $[0, t_1]$. In Section 2 we have solved this problem for the second order discrete-time systems by defining the sequence of measurements at discrete-time instants $t = 0, 1, 2, ..., 2n-1$. By analogy way in the second order continuous-time system (11) by computing higher order derivatives of the continuous-time measurements (12) at a point $t = 0$ we have

$$y(0) = C\,x(0), \ y'(0) = C\,x'(0) = Cx_1, \ y''(0) = CS_0 x_0 + CP_0 x_1 + CM_0 u(0),$$
$$y'''(0) = CS_1 x_0 + CP_1 x_1 + CM_1 u(0) + CM_0 u'(0), \ldots, y^{(2n-1)}(0) = CS_{2n-3} x_0 \quad (13)$$
$$+ CP_{2n-3} x_1 + CM_{2n-3} u(0) + CM_{2n-2} u'(0) + \cdots + CM_0 u^{(2n-3)}(0),$$

where in this case $S_0 = A_0$, $P_0 = A_1$ and denote $S_k = P_{k-1} A_0$; $P_k = S_{k-1} + P_{k-1} A_1$, $k = 1, 2, ...$ ; $M_k = A_0 M_{k-2} + A_1 M_{k-1}$ ($M_0 = B, M_1 = A_1 B$), $k = 2, 3, 4, ...$

By analogy with the second order discrete-time system, (13) can be rewritten in the following matrix form:

$$\begin{bmatrix} y(0) \\ y'(0) \\ y''(0) - CM_0 u(0) \\ \vdots \\ y^{(2n-1)}(0) - CM_{2n-3} u(0) - CM_{2n-2} u'(0) - \cdots - CM_0 u^{(2n-3)}(0) \end{bmatrix}_{(2np) \times 1} = \mathcal{O}(A_0, A_1, C) \cdot \begin{bmatrix} x_0 \\ x_1 \end{bmatrix} \quad (14)$$

where $\mathcal{O} = \mathcal{O}(A_0, A_1, C)$ is the observability matrix already defined in (8). Thus, for determining of the initial values $x_0, x_1$ uniquely from (14) it is necessary and sufficient that the observability matrix has full rank, i.e. rank $\mathcal{O} = 2n$. It would be noted that taking more higher-order derivatives in (14) it is impossible to increase the rank of the observability matrix $\mathcal{O}(A_0, A_1, C)$, because by the Cayley–Hamilton theorem (Kurosh,

A.G. 1972) for $m \geq 2n$ we have $A^m = \sum_{j=0}^{2n-1} \beta_j A^j$ and so the additional equations would be linearly dependent on the previously defined $2n$ equations of (14).

Thus, we have proved the following theorem for the validity of Kalman's type criterion.

**Theorem 3.1** For observability of the linear second order continuous-time system (11) with measurements (12) it is necessary and sufficient that the observability matrix has full rank.

**Remark 3.1** It can be noted out that for a proof of Theorem 3.1 we can use an idea of construction of the discrete-approximation problem for second order continuous system (11); at first we introduce the first and second order difference operators

$$\Delta x(t) := \frac{1}{h}[x(t+h) - x(t)], \ \Delta^2 x(t) := \frac{1}{h}[\Delta x(t+h) - \Delta x(t)], \ t = 0, h, 2h, ...,$$

where $h$ is a step on the $t-$axis and $x(t) \equiv x_h(t)$ is a grid functions on a uniform grid. We define the following discrete-approximation problem associated with the system (11)

$$\Delta^2 x(t) = A_0 x(t) + A_1 \Delta x(t) + Bu(t), \ t = 0, h, ...; \ x(0) = x_0, \ \Delta x(0) = x_1. \quad (15)$$

Note that by using the definition of $s$ th-order difference operator

$$\Delta^s x(t) = \frac{1}{h^s} \sum_{k=0}^{s} (-1)^k C_s^k x(t + (s-k)h), \ C_s^k = \frac{s!}{k!(s-k)!}, \ t = 0, 1, 2, ...$$

we have associated with the linear state-space system (11) the following discrete-approximate equation

$$\Delta^s x(t) = A_0 \Delta^{s-2} x(t) + A_1 \Delta^{s-1} x(t) + B \Delta^{s-2} u(t), \ s = 2, 3, ..., 2n+1.$$

We proceed by analogy with the preceding derivation of Section 2; only in this case in the equations (9) rewritten for second order discrete-approximation-time system (15) instead of $\Delta^s y_0$, $\Delta^s u_0, x_0, x_1, y_0$ we take $\Delta^s y(0), \Delta^s u(0), x(0), \Delta x(0), y(0)$, respectively. Then derivation of the conditions (14) is implemented by passing to the limit here (at least formally), whenever the discrete steps tend to zero.

### 4. Controllability of Second Order Discrete Systems

It should be noted that controllability is an important property of a control system, and besides controllability property plays a crucial role in many control problems, such as stabilization of unstable systems by feedback, or optimal control.

Now, we consider again a linear discrete-time invariant control system defined In Section 2:

$$x_{t+2} = \tilde{A}_0 x_t + \tilde{A}_1 x_{t+1} + \tilde{B} u_t, \ x_0 = a_0, \ x_1 = a_1. \tag{16}$$

Thus, the linear second order discrete-time system modelled by (16) *is controllable* if it is possible to find a control sequence such that the given point can be reached from any initial state in a finite time. In other words, we want to transfer our system from any initial state $x_0 = a_0$, $x_1 = a_1$ to any desired final state $x_{t_1} = x_f$ in a finite time $t_1 < \infty$.

We notice that a dynamic system is guided from an initial state to a desired state by manipulating the control variables. Therefore, when addressing the controllability problem, one needs to consider only the state equation (16), often represented by $(\tilde{A}_0, \tilde{A}_1, \tilde{B})$. The conditions of controllability are imposed on $(\tilde{A}_0, \tilde{A}_1, \tilde{B})$. We continue our techniques used in the study of observability problems of the preceding sections; by repeatedly substitutions here we derive the following set of equations

$$x_2 = S_0 x_0 + P_0 x_1 + M_0 u_0; \; x_3 = S_1 x_0 + P_1 x_1 + M_1 u_0 + M_0 u_1; \; x_4 = S_2 x_0 + P_2 x_1 + M_2 u_0 + M_1 u_1 + M_0 u_2;$$
$$\cdots ; \; x_{n+1} = S_{n-1} x_0 + P_{n-1} x_1 + M_{n-1} u_0 + M_{n-2} u_1 + \cdots + M_0 u_{n-1}. \tag{17}$$

In a general case, as before according to (17) where the input $u_t$ is a vector of dimension $r$, the repetition of the same procedure give us

$$x_{n+1} - S_{n-1} x_0 - P_{n-1} x_1 = \begin{bmatrix} M_0 \vdots M_1 \vdots M_2 \vdots \cdots \vdots M_{n-1} \end{bmatrix} \begin{bmatrix} u_{n-1} \\ \vdots \\ u_1 \\ u_0 \end{bmatrix} \tag{18}$$

where $\begin{bmatrix} M_0 \vdots M_1 \vdots M_2 \vdots \cdots \vdots M_{n-1} \end{bmatrix}$ is a matrix of dimension $n \times (r \cdot n)$. In what follows we denote it by $\mathcal{C}(\tilde{A}_0, \tilde{A}_1, \tilde{B}) = \mathcal{C}^{n \times (n \cdot r)}(\tilde{A}_0, \tilde{A}_1, \tilde{B})$ and define it as the controllability matrix:

$$\mathcal{C}(\tilde{A}_0, \tilde{A}_1, \tilde{B}) = \begin{bmatrix} M_0 \vdots M_1 \vdots M_2 \vdots \cdots \vdots M_{n-1} \end{bmatrix}, \quad \mathcal{C}(\tilde{A}_0, \tilde{A}_1, \tilde{B}) \begin{bmatrix} u_{n-1} \\ \vdots \\ u_1 \\ u_0 \end{bmatrix}_{(nr) \times 1} = x_{n+1} - S_{n-1} x_0 - P_{n-1} x_1 \tag{19}$$

This system consists of $n$ linear algebraic equations in $r \cdot n$ unknowns for $n$ number of $r$-dimensional vector components of $u_0, u_1, ..., u_{n-1}$. Then these equations will have a solutions for any given vector $x_f$ if and only if the matrix has full rank, i.e.

$\text{rank} \mathcal{C}(\tilde{A}_0, \tilde{A}_1, \tilde{B}) = n$.

These results are summarized in the following theorem.

**Theorem 4.1** The linear discrete-time system (16) is controllable if and only if $\text{rank} \mathcal{C}(\tilde{A}_0, \tilde{A}_1, \tilde{B}) = n$, where the controllability matrix is defined by (19).

## 5. Controllability of Second Order Continuous Systems

As is known, the basic concepts of controllability with continuous systems play an essential role in the solutions of many important different optimal control problems. Roughly speaking, as in the discrete systems we remind that controllability of continuous systems generally means that it is possible to steer dynamical system from an arbitrary initial state to an arbitrary final state using the set of admissible controls.

Suppose we have a simplified problem, where the input is a scalar and input matrix $B$ is a vector, usually denoted by $b$:

$$x'' = A_0 x + A_1 x' + bu, \quad x(0) = x_0, x'(0) = x_1. \tag{20}$$

A system is said to be (state) controllable at time $t = 0$, if there exists a finite $t_1 > 0$ such that for any $x_0$ and any $x_1$, there exists an input $u(t), t \in [0, t_1]$ that will transfer the state $(x(0), x'(0))$ to the state $x(t_1)$ at time $t_1$, otherwise the system is said to be uncontrollable at time $t = 0$ (we note that for a fist order linear continuous systems with an arbitrary initial states this notion is also known as complete controllability). Suppose that $x(\cdot)$ and $u(\cdot)$ are continuously differentiable functions of orders $n+1$ and $n-1$ respectively. By sequentially differentiation of the equation (20) we derive the following system of equations

$$x'' = A_0 x + A_1 x' + bu; \quad x''' = A_0 x' + A_1 x'' + bu' = A_1 A_0 x + (A_0 + A_1^2) x' + A_1 bu',$$

$$x^{(IV)} = (A_0^2 + A_1^2 A_0) x + (A_1 A_0 + A_0 A_1 + A_1^3) x' + (A_0 + A_1^2) bu + A_1 bu' + bu'', \ldots.$$

which in terms of designations (4), (5) taking $M_0 = b, M_1 = A_1 b$ can be rewritten in more relevant form

$$x'' = S_0 x + P_0 x' + M_0 u; \quad x''' = S_1 x + P_1 x' + M_1 u + M_0 u', \quad x^{(IV)} = S_2 x + P_2 x' + M_2 u + M_1 u' + M_0 u'',$$
$$\cdots, \quad x^{(n+1)} = S_{n-1} x + P_{n-1} x' + M_{n-1} u + M_n u' + \cdots + M_0 u^{(n-1)}. \tag{21}$$

In fact, in this case the controllability matrix

$$\mathcal{C} = \mathcal{C}(A_0, A_1, b) = [M_0 \vdots M_1 \vdots M_2 \vdots \cdots \vdots M_{n-1}], \quad M_0 = b, M_1 = A_1 b \tag{22}$$

is a square matrix of size $n \times n$. Then if the controllability matrix (22) is non-singular, i.e. its determinant is nonzero, then equations (21) has the unique solution for the input sequence given by

$$\begin{bmatrix} u^{(n-1)}(t) \\ \vdots \\ u'(t) \\ u(t) \end{bmatrix} = \mathcal{C}^{-1} (x^{(n+1)}(t) - S_{n-1} x(t) - P_{n-1} x'(t)). \tag{23}$$

Note that (23) is valid for any $t \in [0, t_f]$ with an arbitrary finite $t_f$. Thus, the nonsingularity of the controllability matrix implies the existence of the scalar input function $u(t)$ and its $n-1$ order derivatives, for any $t < t_f < \infty$. For a vector input system $x'' = A_0 x + A_1 x' + Bu$, where $B$ is $n \times r$ matrix, $u(t) \in \mathbb{R}^r$ is an input vector, the above discussion produces the same relation as (23) with the controllability matrix $\mathcal{C}^{n \times (n \cdot r)}(A_0, A_1, B)$ and with matrices $M_0 = B, M_1 = A_1 B$. In addition, we recall that by Kronecker-Capelli theorem (see, for example (Kurosh, A.G. 1972)) a system of linear equations with controllability matrix $\mathcal{C}^{n \times (n \cdot r)}(A_0, A_1, B)$ is compatible if and only if the rank of the coefficient matrix $\mathcal{C} = \mathcal{C}^{n \times (n \cdot r)}(A_0, A_1, B)$ is equal to that of the augmented matrix $\overline{\mathcal{C}}$ obtained from $\mathcal{C}$ by adding the column of free terms:

$$\operatorname{rank} \mathcal{C}(A_0, A_1, B) = \operatorname{rank}\left[\mathcal{C}(A_0, A_1, B) \vdots d(t)\right], \quad d(t) = x^{(n+1)}(t) - S_{n-1} x(t) - P_{n-1} x'(t).$$

On the other hand a solution of (23) with $\mathcal{C}^{n \times (n \cdot r)}(A_0, A_1, B)$ for any desired state at $t$ and for any $d(t) = x^{(n+1)}(t) - S_{n-1} x(t) - P_{n-1} x'(t)$ exists if and only if $\operatorname{rank} \mathcal{C}(A_0, A_1, B) = n$.

We have obtained the following theorem.

**Theorem 5.1** The linear continuous-time system $(A_0, A_1, B)$ is controllable if and only if the controllability matrix has full rank, i.e. $\operatorname{rank} \mathcal{C}(A_0, A_1, B) = n$.

**Corollary 5.1** For a second order continuous-time linear state-space system $(O_{n \times n}, A_1, B)$ the controllability matrix $\mathcal{C}(O_{n \times n}, A_1, B)$ has the Kalman criterion form $\left[B \vdots A_1 B \vdots A_1^2 B \vdots \cdots \vdots A_1^{n-1} B\right]$.

*Proof.* Recall that $A_0 = O_{n \times n}$ and by the relations (5) $M_k = A_1 M_{k-1}$ $(M_0 = B, M_1 = A_1 B)$, $k = 2, 3, \ldots, n-1$. Then in view of recurrence relation we have

$$M_k = A_1 M_{k-1} = A_1^2 M_{k-2} = \ldots = A_1^k M_0 = A_1^k B. \qquad (24)$$

It means that taking $k = 2, 3, \ldots, n-1$ in (24) we derive $M_2 = A_1^2 B, M_3 = A_1^3 B, \ldots, M_{n-1} = A_1^{n-1} B$. Hence, substituting these relations $(M_0 = B, M_1 = A_1 B)$ in our controllability matrix $\mathcal{C}(O_{n \times n}, A_1, B)$ we have the desired result.

**Corollary 5.2** Suppose that the state and input vector functions of the second order continuous-time linear state-space system $(A_0, O_{n \times n}, B)$ have $2n+1$ and $2n-1$ order derivatives, respectively. Then the controllability matrix $\mathcal{C}(A_0, O_{n \times n}, B)$ has the Kalman criterion form $\left[ B \vdots A_0 B \vdots A_0^2 B \vdots \cdots \vdots A_0^{n-1} B \right]$.

*Proof.* In this case $A_1 = O_{n \times n}$. In terms of matrices of relations (4),(5) we have $P_k = S_{k-1}$ $(S_0 = A_0)$, $M_k = A_0 M_{k-2} (M_0 = B)$. It can be easily checked that for odd and even subindices here $S_{2k} = A_0^{k+1}$,

$P_{2k+1} = A_0^{k+1}$, $M_{2k} = A_0^k B$ and $M_{2k+1} = O_{n \times r}$, $S_{2k+1} = P_{2k+1} = O_{n \times n}$. Then rewriting (21) for the case $M_0 = B$ as a result of differentiation both left and right sides of $x'' = A_0 x + Bu$ we derive that

$$x'' = S_0 x + M_0 u;\ x''' = P_1 x' + M_0 u',\ x^{(IV)} = S_2 x + M_2 u + M_0 u'',\ \ldots,\ x^{(2n)} = S_{2n-2} x + M_{2n-2} u$$
$$+ M_{2n-4} u'' + \cdots + M_0 u^{(2n-2)},\ x^{(2n+1)} = P_{2n-1} x + M_{2n-2} u' + M_{2n-4} u''' + \cdots + M_0 u^{(2n-1)}. \quad (25)$$

Now last two equations of the system (25) means that

$$x^{(2n)} = A_0^n x + A_0^{n-1} B u + A_0^{n-2} B u'' + \cdots + B u^{(2n-2)},\ x^{(2n+1)} = A_0^n x' + A_0^{n-1} B u' + A_0^{n-2} B u''' + \cdots + B u^{(2n-1)}.$$
(26)

Hence from relations (26) we deduce two system of linear equations with unique controllability matrix $\mathcal{C}(A_0, B) = \mathcal{C}^{n \times (n \cdot r)}(A_0, O, B) = \left[ B \vdots A_0 B \vdots A_0^2 B \vdots \cdots \vdots A_0^{n-1} B \right]$

$$\mathcal{C}(A_0, B) \begin{bmatrix} u^{(2n-2)}(t) \\ \vdots \\ u''(t) \\ u(t) \end{bmatrix} = x^{(2n)}(t) - A_0^n x(t);\quad \mathcal{C}(A_0, B) \begin{bmatrix} u^{(2n-1)}(t) \\ \vdots \\ u'''(t) \\ u'(t) \end{bmatrix} = x^{(2n+1)}(t) - A_0^n x'(t)$$

**Remark 5.1** We recall that these observability results can be obtained from duality relation between observability and controllability. That is, controllability and observability are dual properties, because the controllability of $(A_0, O_{n \times n}, B)$ is the same as the observability of $(A_0^T, O_{n \times n}, C^T)$. Consequently, the relation $x'' = A_0 x + B u$, $y = Cx$ is observable if and only if the dual system $x'' = A_0^T x + C^T u$ with output vector $y = B^T x$ is controllable, where a prime $T$ denotes the adjoint transformation. But according to the Theorem 4.1 for controllability of the last dual system the controllability $n \times (np)$ matrix $\left[ C^T \vdots A_0^T C^T \vdots (A_0^2)^T C^T \vdots (A_0^3)^T C^T \vdots \cdots \vdots (A^{n-1})^T C^T \right]$ should be have full rank. We recall that for a first order dual systems $A_0^T$ should be replaced by $-A_0^T$. In general, it can be easily checked that for the system $(A_0, A_1, B)$ the dual system is

$$x'' = A_0^T x - A_1^T x' + C^T u.$$

**Example 5.1**. Suppose we have the linear continuous-time system

$$x'' = \begin{bmatrix} 1 & 0 & 2 \\ 2 & 1 & -1 \\ 3 & 0 & -2 \end{bmatrix} x + \begin{bmatrix} 0 & 3 & 1 \\ 4 & 2 & 1 \\ 1 & -2 & 0 \end{bmatrix} x' + \begin{bmatrix} 1 \\ 0 \\ 2 \end{bmatrix} u, \quad A_0 = \begin{bmatrix} 1 & 0 & 2 \\ 2 & 1 & -1 \\ 3 & 0 & -2 \end{bmatrix}, A_1 = \begin{bmatrix} 0 & 3 & 1 \\ 4 & 2 & 1 \\ 1 & -2 & 0 \end{bmatrix}, B = b = \begin{bmatrix} 1 \\ 0 \\ 2 \end{bmatrix}$$

For construction of the controllability matrix $\mathcal{C}(A_0, A_1, b) = [M_0 \vdots M_1 \vdots M_2]$ we should compute the matrices $M_1, M_2$. As is known, $M_1 = A_1 M_0$, $M_2 = A_0 M_0 + A_1 M_1$, $M_0 = b$. Then

$$M_0 = \begin{bmatrix} 1 \\ 0 \\ 2 \end{bmatrix}, M_1 = \begin{bmatrix} 2 \\ 6 \\ 1 \end{bmatrix}, M_2 = \begin{bmatrix} 24 \\ 21 \\ -11 \end{bmatrix} \text{ and } \mathcal{C}(A_0, A_1, b) = [M_0 \vdots M_1 \vdots M_2] = \begin{bmatrix} 1 & 2 & 24 \\ 0 & 6 & 21 \\ 2 & 1 & -11 \end{bmatrix}.$$

Since $\det \mathcal{C}(A_0, A_1, b) = -291 \neq 0$ we can conclude that $\text{rank} \mathcal{C}(A_0, A_1, b) = 3$ i.e. $\text{rank} \mathcal{C}(A_0, A_1, b) = 3 = n$ implies that the system under consideration is controllable.

**Lemma 5.1** The transfer function of the second order continuous-time linear state-space system $x'' = A_0 x + Bu$ with the corresponding measurements $y = Cx$ has the form

$$H(s) = C \left[ s^2 E - A_0 \right]^{-1} B.$$

*Proof.* Introducing a new variable $v = \begin{bmatrix} x \\ x' \end{bmatrix}$ it is easy to see that our system can be converted to the first order system $v' = \tilde{A}v + \tilde{B}u,$ where $\tilde{A} = \begin{bmatrix} O_{n \times n} & E_n \\ A_0 & O_{n \times n} \end{bmatrix}, \tilde{B} = \begin{bmatrix} O_{n \times r} \\ B \end{bmatrix}.$

By analogy $y = \begin{bmatrix} C & O_{p \times n} \end{bmatrix} v$. Then for a new problem the transfer function is

$H(s) = \begin{bmatrix} C & O_{p \times n} \end{bmatrix} \begin{bmatrix} sE_{2n} - \tilde{A} \end{bmatrix}^{-1} \tilde{B}.$ ($E_k$ is $k \times k$ identity matrix). Now it is not hard to compute the inverse matrix of the invertible matrix $\begin{bmatrix} sE_{2n} - \tilde{A} \end{bmatrix}$ partitioned into a block form:

$$\begin{bmatrix} sE_{2n} - \tilde{A} \end{bmatrix}^{-1} = \begin{bmatrix} sE_n & -E_n \\ -A_0 & sE_n \end{bmatrix}^{-1} = \begin{bmatrix} \begin{bmatrix} sE_n - \frac{1}{s}A_0 \end{bmatrix}^{-1} & \begin{bmatrix} s^2E_n - A_0 \end{bmatrix}^{-1} \\ \begin{bmatrix} s^2E_n - A_0 \end{bmatrix}^{-1} & \begin{bmatrix} sE_n - \frac{1}{s}A_0 \end{bmatrix}^{-1} \end{bmatrix}. \tag{27}$$

Then substitution (27) into the transfer function we have immediately

$$H(s) = C \begin{bmatrix} s^2 E - A_0 \end{bmatrix}^{-1} B, \ E = E_n.$$

**Example 5.2** Consider the following system with scalar input and measurements

$$x'' = \begin{bmatrix} 2 & 1 \\ 3 & 4 \end{bmatrix} x + \begin{bmatrix} 1 \\ 2 \end{bmatrix} u; \ y = \begin{bmatrix} 1 & 3 \end{bmatrix} x; \ A_0 = \begin{bmatrix} 2 & 1 \\ 3 & 4 \end{bmatrix}; \ b = \begin{bmatrix} 1 \\ 2 \end{bmatrix}; \ C = \begin{bmatrix} 1 & 3 \end{bmatrix}.$$

According to Lemma 5.1 let $H(s) = C \begin{bmatrix} s^2 E - A_0 \end{bmatrix}^{-1} b$ be the transfer function of this single-input single-output system, where $E$ is $2 \times 2$ identity matrix. Then

$$H(s) = \begin{bmatrix} 1 & 3 \end{bmatrix} \begin{bmatrix} \dfrac{s^2-4}{(s^2-2)(s^2-4)-3} & \dfrac{1}{(s^2-2)(s^2-4)-3} \\ \dfrac{3}{(s^2-2)(s^2-4)-3} & \dfrac{s^2-2}{(s^2-2)(s^2-4)-3} \end{bmatrix} \begin{bmatrix} 1 \\ 2 \end{bmatrix} = \dfrac{7s^2-5}{s^4-6s^2+5}$$

and consequently, there is no zero-pole cancellations in the transfer function of our a single-input single-output system, that is the system is both controllable and observable (i.e., the poles are $\pm\sqrt{5}$; $\pm 1$ and zero is $\pm\sqrt{5/7}$). Notice that systems that satisfy this relationship are called proper. To get the complete answer we have to go to a state space form and examine the controllability and observability matrices;

$\text{rank}\,\mathcal{O}(A_0, O_{n\times n}, C) = \text{rank}\begin{bmatrix} 1 & 3 \\ 11 & 13 \end{bmatrix} = n = 2$ and the system is observable. On the other hand, $\mathcal{C}(A_0, O_{n\times n}, b) = \begin{bmatrix} 1 & 4 \\ 2 & 11 \end{bmatrix}$, $\det \mathcal{C}(A_0, O_{n\times n}, b) = 3 \neq 0$ and our system is controllable.

## 6. Conclusion

In this paper, the issue on the observability and controllability criteria for a class of Second Order Linear Time Invariant Systems has been addressed. Several sufficient and necessary conditions for observability and controllability of such systems have been established in the form of Kalman's type conditions. By sequentially differentiation of state and output vector-functions with scalar and vector inputs deriving different systems of linear algebraic equations are also discussed, transfer function is constructed. It is worth noticing that this note considers the controllability and observability for a class of linear time invariant systems. Many issues are still untouched on more general linear systems, for example, the higher order linear systems with variable coefficients. Numerical examples are presented to illustrate the theoretical results and to show the effectiveness of the proposed results.